\documentstyle[12pt]{article}

  \newcommand{\const}{\rm const}
  \newcommand{\Var}{\rm Var}
  \newcommand{\supp}{\rm supp}

  \newcommand{\Dom}{\rm  Dom}

   \newcommand{\Sub}{\rm Sub}

\begin{document}

   \begin{center}

{\bf Averaging of random variables and fields. }\\

\vspace{5mm}

{\bf M.R.Formica, \ E.Ostrovsky, \ L.Sirota. }

 \end{center}

 \vspace{5mm}

\ Universit\`{a} degli Studi di Napoli Parthenope, via Generale Parisi 13, Palazzo Pacanowsky, 80132,
Napoli, Italy. \\
e-mail: \ mara.formica@uniparthenope.it \\

 \ Department of Mathematics and Statistics, Bar-Ilan University,\\
59200, Ramat Gan, Israel. \\
e-mail: \ eugostrovsky@list.ru\\

 \ Department of Mathematics and Statistics, Bar-Ilan University,\\
59200, Ramat Gan, Israel. \\
e-mail: \ sirota3@bezeqint.net \\

\vspace{5mm}

  \begin{center}

   {\bf Abstract}

  \end{center}

\vspace{4mm}

 \ We will prove that by averaging of random variables (r.v.) and random fields (r.f.) its tails of distributions do not
 increase in comparison with  the tails of source variables,  essentially or almost exact,  under  very weak conditions. \par

\vspace{5mm}

 \hspace{3mm} {\it  Key words and phrases.}  Random variables (r.v.) and random fields (r.f.), spatial and conditional averaging, probability,
 Lebesgue - Riesz, Yudovich,  Grand Lebesgue norm and spaces, Euclidean space, martingales, uniform integrability,
 tail of distribution,  measurable functions,
Young - Fenchel (Legendre) transform,  Young's inequality,  slowly varying function, natural function,  normalized sums, filtration,
 generating function, upper and lover estimation,  conditional expectation, rearrangement invariant (r.i.) space, examples,
 square variation.\par

\vspace{5mm}

\section{Preliminary. Definitions. Notations.}

\vspace{5mm}

 \hspace{3mm}  Let $ \ (\Omega = \{\omega\}, \cal{B}, {\bf P}) \ $ be non - trivial probability space with expectation $ \ {\bf E} \ $
 and Variation $ \ {\bf \Var} \ $  and let $ \ (X = \{x\}, M, \mu)  \ $ be another probability space: $ \ \mu(X) = 1.  \ $  Let also
 $ \  \xi = \xi(x) = \xi(x,\omega), \ x \in X, \ \omega \in \Omega \ $ be numerical valued separable bi - measurable random field. \par
 \ There are at last two version to introduce the notion of  an {\it averaging} in the probability theory. \par

 \vspace{3mm}

 \hspace{3mm} {\bf A.} Define

 \vspace{3mm}

\begin{equation}  \label{first version}
\xi[X] \stackrel{def}{=} \int_X \xi(x) \ \mu(dx) \ -
\end{equation}

\vspace{3mm}

the spatial averaging.\par

\vspace{4mm}

\hspace{3mm}{\bf B. \ Conditional averaging.}  Let $ \  F  \  $ be some sigma - subalgebra of the source one $ \  \cal{B}.  \ $ Define for
arbitrary integrable numerical valued r.v.  $ \ \eta \ $ the ordinary conditional expectation

\begin{equation}  \label{first version}
\eta\{F\} \stackrel{def}{=} {\bf E}\eta/F.
\end{equation}

\vspace{5mm}

 \hspace{3mm}  {\bf Our target in this preprint is to deduce the exact norm and tail of distribution estimations for both this averaging  in the terms
 of tails and norms of the source datum.  \par

 \ We state as it is announced  that by both the introduced  averaging: of random variables (r.v.) and random fields (r.f.) its tails of distributions
 do not essentially exact increase in comparison with  the tails of source variables,  under  very weak conditions. } \par

 \vspace{5mm}

 \ Recall that the so - called {\it tail function}  $ \ T_{\tau}(t) \ $  for the  numerical valued r.v. $ \ \tau \ $ is defined as follows

\begin{equation} \label{tail fun}
\ T_{\tau}(t) \stackrel{def}{=} {\bf P} (|\tau| \ge t), \ t \ge 0.
\end{equation}

 \ Correspondingly, the tail function for the numerical valued random field $ \ \zeta = \zeta(x) = \zeta(x,\omega) \ $ is understood in the uniform sense

$$
T_{\zeta}(t) \stackrel{def}{=} \sup_{x \in X} T_{\zeta(x)}(t), \ t \ge 0,
$$
is of course the last $ \ \sup \ $  is not trivial. \par

\vspace{5mm}

\section{Averaging of random fields.}

\vspace{5mm}

 \hspace{3mm} {\bf Proposition 2.1.}  \ Norm estimation. \par

 \vspace{3mm}

  \ Let $ \ (G, \ ||\cdot||G)  \ $ be arbitrary Banach functional  space defined on the random variables
  (measurable functions)  $ \ \zeta: \ \Omega \to R, \ $  such that

$$
\forall x \in  X \  \Rightarrow   \xi(x) \in G.
$$
 \ Suppose also that the function $ \ \beta(x) := ||\xi(x)||G, \ x \in X \ $ is measurable. Then by virtue of integral  version of
 triangle inequality

$$
|| \ \xi[X] \ ||G = ||\int_X \xi(x) \ \mu(dx)||G  \le
$$

\begin{equation}  \label{first estimate}
\int_X ||\xi(x)||G \ \mu(dx)  =  \int_X \beta(x) \ \mu(dx),
\end{equation}
 if of course the last integral is finite.\par

 \vspace{4mm}

 \ For instance, denote by $ \ L_p(\Omega,P) = L_p(\Omega), \ p \ge 1  \ $  the classical Lebesgue -Riesz space consisting  of all the
 numerical valued r.v.  $   \ \nu  \  $  having a finite norm

$$
||\nu||_{\Omega,p} = ||\nu||L_p(P) \stackrel{def}{=}  \left[ \ {\bf E} \ |\nu|^p \ \right]^{1/p} =
 \left[ \ \int_{\Omega} | \nu(\omega)|^p \ \ {\bf P}(d\omega)  \ \right]^{1/p}.
$$

 \ Then under formulated restrictions

\begin{equation}  \label{Lp estimate}
|| \ \xi[X] \ ||L_p(\Omega) \le \int_X ||\xi(x)||L_p(\Omega) \cdot \mu(dx).
\end{equation}

\vspace{4mm}

 \   As other spaces  may be considered, e.g. Marcinkiewicz, Lorentz ones. Let us bring an interest (in our opinion) example.  The
 particular case of the  classical Lorentz  quasi - norm (and correspondent space) for the random variable $ \  \phi  \ $  is defined
 as follows

$$
||\phi|| L^{q,\infty} \stackrel{def}{=} \sup_{t \ge 0} \left[ \  t^q \ T_{\phi}(t) \ \right]^{1/q}.
$$
 \ To be more precisely, there exists an actually norm which is equivalent  to the described before  quasi - norm \ \cite{Stein}. \par

  \ We derive applying these spaces. \par

\vspace{4mm}

 {\bf  Proposition 2.2.} Suppose for the introduced r.f $ \ \xi(x) \ $

$$
\sup_{x \in X} T_{\xi(x)}(t) \le t^{-q}, \ t \ge 1, \ q = \const > 1.
$$

 \ Then

$$
\exists C = C(q) < \infty, \ \Rightarrow \ T_{\xi[X]}(t) \le  C \ t^{-q}, \ t \ge 1.
$$

\vspace{4mm}

 \ Let us consider an important case for the r.v. having in general case an exponential decreasing tail of distribution.\par

\begin{center}

 \  {\sc  Grand Lebesgue Spaces.} \par

\end{center}

\vspace{4mm}

\hspace{3mm} Let $ \  b = \const, 1 < b \le \infty.  \  $ Let also $ \ \psi = \psi(p), \ p \in (1,b) \ $  be certain numerical valued
 measurable strictly positive: $ \  \inf_{p \in [1,b)} \psi(p) > 0   \ $  function, not necessary to be bounded. Denotation:
 $ \ \Dom (\psi) \stackrel{def}{=} \{ \ p: \ \psi(p) < \infty \ \}, \ $

 $$
  (1,b) := \supp (\psi); \  \Psi[1,b) := \{ \  \psi: \ \supp (\psi) = (1,b)  \ \},
 $$

$$
\Psi \stackrel{def}{=} \cup_{b > 1} \Psi[1,b).
$$

\vspace{3mm}

 \hspace{3mm} {\bf Definition 2.1.}, see e.g.  \cite{KosOs}, \cite{Ermakov}, \cite{Fiorenza-Formica-Gogatishvili-DEA2018}.
 \ Let the function $ \ \psi(\cdot) \in \Psi[1,b), \ $  which is named as {\it generating function} for introduced after space.
 \ The  so - called {\it Grand Lebesgue Space} $ \  G \psi \ $
 is defined as a set of all random variables (measurable functions)  $ \ \tau \ $ having a finite norm

\vspace{3mm}

\begin{equation} \label{def GLS norm}
||\tau||G\psi  \stackrel{def}{=}  \sup_{p \in (1,b)} \left\{ \ \frac{||\tau||L_p(\Omega)}{\psi(p)} \ \right\}.
\end{equation}

\vspace{4mm}

 \ The particular case of these spaces  and under  some additional restrictions on the generating function $ \ \psi = \psi(p) \ $
  correspondent to the so - called {\it Yudovich spaces,}  see  \cite{Yudovich 1}, \cite{Yudovich 2}.
These spaces was applied at first in the theory of Partial Differential Equations (PDE),  see  \cite{Chen},  \cite{Crippa}. \par

\vspace{3mm}

 \  These spaces are  complete Banach functional rearrangement invariant; they are
 investigated in many works, see e.g. \cite{ErOs}, \cite{Ermakov}, \cite{Fiorenza2}, \cite{Fiorenza-Formica-Gogatishvili-DEA2018},
\cite{fioforgogakoparakoNAtoappear}, \cite{fioformicarakodie2017}, \cite{formicagiovamjom2015}, \cite{Formica Ostrovsky Sirota weak dep},
\cite{KosOs}, \cite{KozOsSir2017}, \cite{KosOs equivalence}, \cite{Liflyand}, \cite{Ostrovsky1}. It is important for us in particular to note
that there is exact of course up to finite multiplicative constant
interrelations  under certain natural conditions on the generating function between belonging the r.v. $ \ \tau \ $ to this
space and it tail behavior. Indeed, assume for definiteness that $ \  \tau \in G \psi  \ $ and moreover $ \ ||\tau||G\psi = 1;  \ $ then

\begin{equation} \label{tail estim}
T_{\tau}(t) \le  \exp \{ \  -  h^*(\ln t)  \  \}, \ t \ge e,
\end{equation}
where $ \ h(p) = h[\psi](p) := p \ln \psi(p) \ $  and $ \  h^*(\cdot)   \ $  is famous Young - Fenchel (Legendre) transform of the function $ \ h(\cdot): \ $

$$
h^*(u) \stackrel{def}{=} \sup_{p \in \Dom( \psi)} (pu - h(p)).
$$

 \ Inversely, let the tail function $ \ T_{\tau}(t), \ t \ge 0  \ $  be given.  Introduce the following so - called {\it natural function}
 generated by $ \ \tau \ $

\begin{equation} \label{natural function}
 \psi_{\tau}(p) \stackrel{def}{=} \left[ \   p \int_0^{\infty}  \ t^{p-1} \ T_{\tau}(t)  \ dt \ \right]^{1/p} = ||\tau||L_p(\Omega),
\end{equation}
 if it is finite for some value $ \ b \in (1,\infty], \ $ following, it is finite at last for all the values $ \ p \in [1,b). \ $ \par

 As long as

$$
{\bf E} |\tau|^p = p \int_0^{\infty}  \ t^{p-1} \ T_{\tau}(t)  \ dt \ = \psi^p_{\tau}(p), \ p \in [1,b),
$$
 we conclude that if the last {\it natural} for the r.v. $ \ \tau \ $ function $ \ \psi_{\tau}(p) \ $ is finite inside some non - trivial
 segment $ \ p \in [1,b), \ 1 < b \le \infty, \ $   then

$$
\tau \in G\psi_{\tau}; \ \hspace{3mm} ||\tau||G\psi_{\tau} = 1.
$$

\vspace{5mm}

 \ Further, let the estimate (\ref{tail estim}) be given. Suppose in addition   that the generating function $ \ \psi = \psi(p), \ p \in \Dom(\psi) \ $
 is continuous and suppose in the case when $ \ b = \infty \ $

\begin{equation} \label{lim 0}
\lim_{p \to \infty} \frac{\psi(p)}{p} = 0.
\end{equation}

 \ Then the r.v. $ \ \tau \ $ belongs to the Grand Lebesgue Space $ \ G \psi: \    \ $

\begin{equation} \label{inverse est}
||\tau||G\psi \le K[\psi] < \infty,
\end{equation}
see e.g.   \cite{KosOs equivalence}.\par

 \ These conditions on the generating function $ \ \psi(\cdot) \ $  are satisfied for example for the functions  $ \ \psi_{m,L}(p) \ $ of the form

\begin{equation} \label{psi mL}
\psi_{m,L}(p) \stackrel{def}{=} p^{1/m} \ L(p), \ m = \const > 1, \ b = \infty,
\end{equation}
 where $ \ L = L(p) \ $ be some continuous strictly positive {\it  slowly varying }  at infinity function such that

\begin{equation} \label{theta condition}
\forall \theta > 0 \ \Rightarrow \sup_{p \ge 1} \left[  \ \frac{L(p^{\theta})}{L(p)} \ \right] = C(\theta) < \infty.
\end{equation}

 \ For instance, $ \ L(p) = [\ln(p+1) ]^r, \ r \in R. \ $ \par

\vspace{3mm}

 \ We conclude that under formulated restrictions the r.v. $ \ \tau \ $ belongs to the space $ \ G\psi_{m,L}: \ $

\begin{equation} \label{Gpsi m L}
\sup_{p \ge 1}  \left\{ \ \frac{||\tau||_{p,\Omega}}{\psi_{m,L}(p)}  \ \right\} = C(m,L) < \infty
\end{equation}

\ if and only if

\begin{equation} \label{tail behav}
T_{\tau}(u) \le \exp \left( \  - C_2(m,L) \ u^m /L(u) \ \right), \ u \ge e, \  \exists \ C_2(m,L) > 0.
\end{equation}

\vspace{4mm}

 \ A very popular example of these spaces  forms the so - called subgaussian space  $ \ \Sub = \Sub(\Omega); \ $ it consists on the
 subgaussian random variables, for which $ \  \psi(p) = \psi_2(p):= \sqrt{p}: \ $

 \begin{equation} \label{def sub}
 ||\tau|| \Sub = ||\tau||G\psi_2 \stackrel{def}{=} \sup_{p \ge 1} \left[   \  \frac{||\tau||_{p, \Omega}}{\sqrt{p}}  \ \right].
 \end{equation}
\ The r.v. $ \ \tau \ $ belongs to the subgaussian space $ \ \Sub(\Omega) \ $ iff

\begin{equation} \label{tail sub}
\exists C > 0 \ \Rightarrow \  T_{\tau}(u)  \le \exp(- C u^2), \ u \ge 0.
\end{equation}

\vspace{3mm}

 \ {\bf Remark 2.1.} As a rule, on the the r.v.  $ \ \tau \ $ from the spaces $ \ G\psi_{m,L} \ $ is imposed the condition
 of {\it centering:} $ \ {\bf E} \tau = 0. \ $ \par

 \vspace{4mm}

 \ Another examples. Suppose that the r.v. $ \ \tau \ $ be such that

$$
T_{\tau}(t) \le T^{\beta,\gamma,L}(t), \ \beta > 1, \ \gamma > -1, \ L = L(t),
$$
where

$$
T^{\beta,\gamma,L}(t) \stackrel{def}{=} t^{-\beta} \ (\ln t)^{\gamma} \ L(\ln t), \ t \ge e
$$
and $ \ L = L(t), \ t \ge e \ $ be as before slowly varying at infinity positive continuous function.
It is known  \ \cite{KosOs equivalence} \ that as $ \ p \in [1,\beta) \ $

\begin{equation} \label{beta gamma}
\psi_{\tau}(p) = ||\tau||_p \le C_1(\beta,\gamma,L) \ (\beta - p)^{-(\gamma + 1)/\beta} \ L^{1/\beta}(1/(\beta - p)),
\end{equation}
and conversely, if the relation (\ref{beta gamma}) there holds, then

$$
T_{\tau}(t) \le  C_7(\beta,\gamma, L) \ T^{\beta,\gamma+1,L}(t).
$$

 \ Herewith both this estimations are unimprovable. \par

\vspace{5mm}

 \ {\bf Theorem 2.1.} Denote for the random field $ \ \xi(x), \ x \in X \ $

$$
Q(t) := \sup_{x \in X} T_{\xi(x)}(t), \ t \ge 0,
$$

$$
g_0(p) := \left[ \ p \ \int_0^{\infty} \ t^{p-1} \ Q(t) \ dt \ \right]^{1/p},
$$

$$
g(p) := p \ln g_0(p),
$$
 and suppose $ \ \exists b \in (1,\infty] \ \Rightarrow g(b) < \infty. \ $ We assert

\vspace{4mm}
\begin{equation} \label{ran field est}
T_{\xi[X]}(t) \le \exp \left( \ - g^*(\ln t) \ \right), \ t \ge e.
\end{equation}

 \vspace{4mm}

 \ {\bf Proof.} We have

$$
T_{\xi(x)} (t) \le Q(t), \ x \in X, \ t \ge 0.
$$

 \ Therefore the random variables $ \ \xi(x), \ x \in X \ $ belongs to the  unit ball of the Grand Lebesgue Space $ \ G  \ g_0: \ $

$$
\sup_{x \in X} ||\xi(x)||G\kappa \le 1.
$$

 \ It remains to apply the proposition (\ref{tail estim}). \par

\vspace{4mm}

 \ {\bf Remark 2.2.}  A more general version: assume that there exists a finite a.e.  and non - negative measurable
 numerical valued function  $ \ \theta = \theta(x), \ x \in X \ $ such that

$$
T_{\xi(x)}(t) \le \exp \left\{ \  - h^*_{\psi}[\ln (t/\theta(x))] \ \right\}, \ t \ge e,
$$

then

$$
||\xi(x)|| G\psi \le C_4[\psi] \ \theta(x), \ C_4[\psi] < \infty,
$$

therefore

$$
||\xi[X]||G \psi \le C_4[\psi] \int_X \theta(x) \ \mu(dx),
$$
with correspondent tail estimate

\begin{equation} \label{integral form}
T_{\xi[X]} \le \exp \left[ \ - h^*_{\psi}(\ln t/C_5)  \ \right], \ t \ge e C_5,
\end{equation}

where

$$
C_5 = C_4[\psi] \ \int_X \theta(x) \ \mu(dx),
$$
if of course the last integral is finite. \par

\hspace{4mm}

 \ {\bf Example 2.1.} Suppose

\vspace{3mm}

\begin{equation} \label{tail behav}
\sup_{x \in X} T_{\xi(x)}(u) \le \exp \left( \  - \ u^m /L(u) \ \right), \ u \ge e,
\end{equation}

\vspace{3mm}
where  as before $ \  m = \const > 1, \ L = L(p) \ $ is some continuous strictly positive {\it  slowly varying }  at
 infinity function such that

\begin{equation} \label{theta again condition}
\forall \theta > 0 \ \Rightarrow \sup_{p \ge 1} \left[  \ \frac{L(p^{\theta})}{L(p)} \ \right] = C(\theta) < \infty.
\end{equation}

\vspace{3mm}

 \ We conclude that under formulated  above restrictions that $ \ \exists C_5(m,L) > 0 \ \Rightarrow  \ $

\begin{equation} \label{tail exampl}
T_{\xi[X]}(u) \le \exp \left \{ \  - C_5(m,L) \ u^m /L(u) \ \right\}, \ u \ge e.
\end{equation}

 \vspace{3mm}

 \ Obviously, the last estimate (\ref{tail exampl}) is non - improvable, of course. up to multiplicative constant
 $ \ C_5(m,L). \ $ \par

\vspace{4mm}

 \ {\bf Example 2.2.} Suppose now in the introduced notations and restrictions

\vspace{3mm}

\begin{equation} \label{tail  beta}
\sup_{x \in X} T_{\xi(x)}(u) \le T^{\beta,\gamma, L}(t), \ t \ge 1,
\end{equation}

\ We  derive arguing similarly to the previous example

\begin{equation} \label{beta gamma est}
 T_{\xi[X]}(u) \le T^{\beta,\gamma +1, L}( C_6(\beta,\gamma,L) \ t), \ t \ge e, \ C_6 > 0,
\end{equation}
and this estimate is also non - improvable still in the case when the set $ \ X \ $ consists in a single point \ \cite{Liflyand}.\par

\vspace{4mm}

 \ {\bf Example 2.3.} Let us show that the  proposition of theorem 2.1 is not true if we replace  the averaging onto the classical
 {\it normalized} sums of the centered random variables. Let $ \ \{\zeta_i\}, \ i = 1,2,\ldots; \ \zeta = \zeta_1 \ $ be a sequence if independent
 symmetrical distributed r.v. such that

$$
\forall y \ge 0 \ \Rightarrow  \ T_{\zeta}(y) = \exp \left( \ - y^q \ \right), \ q = \const > 2.
$$
 \ Denote

$$
S_n := n^{-1/2} \sum_{i=1}^n\zeta_i; \ \sigma^2 = \sigma^2(q)  \stackrel{def}{=} \Var (\zeta) \in (0,\infty).
$$

 \ We derive by virtue of CLT for all the positive values $ \ y, \ $ say for $ \ y \ge 1 \ $

$$
\sup_n {\bf P} (S_n > y) \ge \lim_{n \to \infty} {\bf P} (S_n > y) \ge \exp \left(-C_0 \ y^2 \right),
$$
 $ \ C_0 = C_0(q) \in (0,\infty); \ $
 so that the  tails of distribution of the r.v. $ \ S_n \ $ are "much heavier" as ones for the source variable $ \ \zeta. \ $ \par

\vspace{5mm}

\section{Conditional averaging.}

\vspace{5mm}

  \hspace{3mm} Let as in the first section (\ref{first version}) $ \  \nu \stackrel{def}{=} \eta\{F\} \stackrel{def}{=} {\bf E}\eta/F, \ $ where
  $ \ F \ $ is some sub - sigma field of the source  sigma algebra  {\bf \cal{B}}. \par

  \vspace{4mm}

 \ Assume that

$$
T_{\eta}(t) \le R(t), \ t \ge 0;
$$
where $ \ R = R(t) \ $ is non - negative bounded: $ \ R(t) \le 1 \ $ non - increasing measurable function such that $ \ R(\infty) = 0.\ $
Introduce the $ \ \Psi \ $ function

$$
\delta(p) = \delta[R](p)  \stackrel{def}{=} \left[ \ p \ \int_0^{\infty} t^{p-1} \ R(t) \ dt \  \right]^{1/p},
$$
and suppose $ \ \delta(\cdot) \in \Psi;  \ $ then

$$
{\bf E} |\eta|^p \le \delta^p(p).
$$

 \ Since the function $ \ g(y) = |y|^p, \ p \ge 1 $ is convex, one can apply the Jensen's inequality for
 the conditional expectations

$$
{\bf E} |\nu|^p \le {\bf E} |\eta|^p \le  \delta^p(p), \ p \in [1,b);
$$
following $ \  \nu \in G \delta  \ $ and we obtain the inequality

\begin{equation} \label{tail v condit}
T_{\nu}(t) \le \exp \left( \ - v^*(\ln  \ t)  \ \right), \ t \ge e,
\end{equation}
where

$$
v(p) = p \ \ln \delta(p), \ p \in \Dom (\delta).
$$

\vspace{4mm}

 \ To summarize. \par

  \ {\bf Proposition 3.1.} We obtain under formulated notations and restrictions

\begin{equation} \label{Gdelta}
||\nu||G\delta  \le 1,
\end{equation}
and hence the estimation (\ref{tail v condit}) there holds.\par

\vspace{4mm}

 \ The examples may be considered alike ones in the previous section, under at the same restrictions. \par

\vspace{3mm}

{\bf Example  3.1.} If

\begin{equation} \label{condit tail}
T_{\eta}(u) \le \exp \left \{ \  - u^m /L(u) \ \right\}, \ u \ge e,
\end{equation}

then

\begin{equation} \label{tail result cond}
T_{\nu}(u) \le \exp \left \{ \  - C_7(m,L) \ u^m /L(u) \ \right\}, \ u \ge e, \ C_7 > 0.
\end{equation}

 \vspace{3mm}

\ {\bf Example 3.2.} Suppose

\vspace{3mm}

\begin{equation} \label{tail condit}
T_{\eta}(u) \le T^{\beta,\gamma, L}(t), \ t \ge 1, \ \beta > 1, \ \gamma > - 1.
\end{equation}

\ We  derive

\begin{equation} \label{beta gamma est}
 T_{\nu}(u) \le T^{\beta,\gamma +1, L}( C_8(\beta,\gamma,L) \ t), \ t \ge e, \ C_8 > 0,
\end{equation}
and this estimate is also essentially non - improvable, see e.g.  \ \cite{Liflyand}.\par

\vspace{3mm}

 \ {\bf Example 3.3.} Let $ \ \Omega = (0,1) \ $ equipped with the classical Lebesgue measure. Define
the r.v. $ \ \eta  = \omega^{-\alpha},  \ \omega  \in (0,1); \ \alpha = \const \in (0,1) \ $  and put

$$
F := \{ \ \emptyset, \ (0,1/2], (1/2,1), (0,1) \ \}.
$$
 \ We have

$$
{\bf E} |\eta|^p = \frac{1}{1 - \alpha p}, \ 1 \le p < 1/\alpha
$$
and $ \ {\bf E} |\eta|^p = \infty, \ p \ge 1/\alpha. \ $ The r.v. $ \ \nu = {\bf E} \eta/F \ $ has a form

$$
\nu(\omega) = 2^{\alpha}/(1 - \alpha), \ \omega \in (0,1/2]
$$
and

$$
\nu(\omega) = (2 - 2^{\alpha})/(1 - \alpha), \ \omega \in (1/2,1).
$$

\vspace{3mm}

 \ Thus, the r.v. $ \ \nu \ $ is {\it bounded}, despite that the source r.v. $ \ \eta \ $ has not all the moments
 $ \ {\bf E} |\eta|^p, \ p \ ge 1/\alpha. \ $ \par

\vspace{5mm}

\section{Applications to the martingale theory.}

\vspace{5mm}

 \begin{center}

 {\bf A.} {\sc Uniform integrable martingales.} \par

 \end{center}

\vspace{4mm}

 \ Suppose that there is  certain {\it filtration} on the source   probability sigma - field $ \ {\bf \cal{B}}, \ $
 i.e. an increasing sequence of sigma - subfields $ \ \{F_n \} \ $
 of canonical one $ \ {\bf \cal{B}}, \ $  and let $ \ (\kappa_n, F_n), \ n = 1,2,\ldots   \ $ be an {\it uniform integrable} martingale:

 $$
 {\bf E} \kappa_m/F_k = \kappa_k, \ 1 \le k \le m
 $$

 \ We conclude  by virtue of theorem J.Doob there exists with probability one a limit

 \begin{equation} \label{lim Doob}
 \kappa \stackrel{def}{=} \lim_{n \to \infty} \kappa_n.
 \end{equation}

 \ Assume that

$$
T_{\kappa}(t) \le  U(t), \ t \ge 0,
$$
 where $ \ U = U(t) \ $ is some non - negative  bounded $ \ U(t) \le 1 \ $ non - increasing function for which
 $ \ \lim U(t) = 0, \ t \to\infty. $
 \  Introduce  as before the following $ \ \Psi - \ $ function

$$
\rho(p) \stackrel{def}{=} \left[ \ p \ \int_0^{\infty} t^{p-1} \ U(t) \ dt \ \right]^{1/p},
$$
and we suppose $ \ \rho \in \Psi(b) \ $ for some value $ \ b \in (1,\infty]. \ $ \par

\vspace{3mm}

 \ We get on the basis of proposition 3.1  using the fact that $ \ \kappa_n = {\bf E} \kappa/F_n  $ \par

\vspace{4mm}

 {\bf  Proposition 4.1.}

\vspace{3mm}

\begin{equation} \label{unif integr}
||\kappa_n||G\rho  \le 1
\end{equation}

\vspace{3mm}

with correspondent tail estimation  (\ref{tail estim}). \par

\vspace{4mm}

 \  Further, define the following finite a.e. random variable

$$
\hat{\kappa} := \sup_n \kappa_n
$$
and a new $ \ \Psi \ - \ $ function

$$
\hat{\rho}(p) := p/(p-1) \cdot \rho(p), \ p \in (1,b).
$$

\vspace{4mm}

 {\bf  Proposition 4.2.}

\vspace{3mm}

\begin{equation} \label{unif sup integr}
|| \hat{\kappa}||G \hat{\rho}  \le 1,
\end{equation}

with correspondent tail estimation. \par

\vspace{4mm}

 {\bf Proof} is quite alike to the previous one; one can apply the famous Doob's inequality

 $$
 {\bf E} |\hat{\kappa}|^p \le \left[ \ \frac{p}{p - 1} \ \right]^p \cdot \rho^p(p)  = \left[ \ \hat{\rho}(p) \ \right]^p.
 $$

 \vspace{4mm}

 \ {\bf Remark 4.1.} The proposition 4.2 allows a simple (particular) inversion. Assume namely that the estimation (\ref{unif sup integr}) is satisfied.
Then the r.v. $ \  \{ \kappa_n \} \ $ are uniform integrable and, following, by  virtue of theorem J.Doob there exists a limit a.e.

\begin{equation} \label{lim ae}
{\bf P} (\lim_{n \to \infty} \kappa_n = \kappa) = 1.
\end{equation}

 \ It follows from (\ref{lim ae}) and   (\ref{unif sup integr}) that in turn

$$
\sup_n ||\kappa_n||G \hat{\rho} < \infty.
$$

\vspace{4mm}

 \begin{center}

 {\bf B.} {\sc General case: a non - uniform integrable martingales.} \par

 \end{center}

\vspace{4mm}

 \ Let again $ \  (\kappa_n, F_n)  \ $ be a martingale; we do not suppose in this subsection that it is
 uniform integrable. Then the sequence $ \ \kappa_n \ $ must be normed; and  as it is noted in \cite{Pekshir Shiryaev}
the natural norming sequence may be choosed as a its {\it square variation}

\begin{equation} \label{square var}
[\kappa]_n^2 \stackrel{def}{=} \kappa_1^2 + \sum_{i=2}^n (\kappa_i - \kappa_{i-1})^2, \ n \ge 2.
\end{equation}
 \ But the square variation sequence is {\it random}; therefore we offer its {\it expectation} as a capacity of
 norming sequence

\begin{equation} \label{sigma norming}
\sigma_n^2 \stackrel{def}{=} {\bf E} [\kappa]_n^2, \ n \ge2; \hspace{3mm} \sigma_1^2 := {\bf E} \kappa_1^2.
\end{equation}

 \ {\it It will be presumed  further} that $ \ \sigma_n \in (0,\infty), \ n \ge 1. \ $ \par

\vspace{3mm}

 \ For instance, if our martingale may be represented as a sum of the centered independent r.v.

$$
\kappa_n = \sum_{j=1}^n \epsilon_j, \hspace{3mm}  F_n  = \sigma\{ \ \epsilon_j, \ 1 \le j \le n  \ \},
$$
such that $ \ \Var (\epsilon_j) < \infty, \ $ then

$$
\sigma^2_n = \sum_{j=1}^n \Var (\epsilon_j).
$$

\vspace{4mm}

 \ Further, suppose as above  that

$$
\sup_n T_{\kappa_n/\sigma_n}(t) \le W(t), \ t \ge 1,
$$
 where $ \  W(\cdot) \ $ is non - negative bounded $ \ W \le 1 \ $  non - increasing function  such that $ \ W(\infty) = 0. \ $
Define

$$
\upsilon(p)  := \left[ \ p \ \int_0^{\infty} t^{p-1} \ W(t) \ dt  \ \right]^{1/p}, \hspace{3mm}
\hat{\upsilon}(p) := p \upsilon(p).
$$

\vspace{4mm}

 \ {\bf Proposition 4.3.} \par

\vspace{4mm}

 \ One can apply the famous  Burkholder's  inequality \ \cite{Burkholder}

$$
||\max_{k = 1,2,\ldots,n} \kappa_k||_p \le p \ \upsilon(p) = \hat{\upsilon}(p),
$$
or equally

$$
||\max_{k = 1,2,\ldots,n} \kappa_k||G \hat{\upsilon} \le 1
$$
with correspondent tail estimate (\ref{tail estim}). \par

 \vspace{5mm}

\section{Concluding remarks.}

\vspace{5mm}

 \hspace{3mm}

 \ It is interest in our opinion to generalize obtained results on the variables, and martingales,
 taking values in Banach spaces. \par

\vspace{6mm}

\vspace{0.5cm} \emph{Acknowledgement.} {\footnotesize The first
author has been partially supported by the Gruppo Nazionale per
l'Analisi Matematica, la Probabilit\`a e le loro Applicazioni
(GNAMPA) of the Istituto Nazionale di Alta Matematica (INdAM) and by
Universit\`a degli Studi di Napoli Parthenope through the project
\lq\lq sostegno alla Ricerca individuale\rq\rq .\par

\end{document}